\documentclass[twoside]{article}
\usepackage{amsmath}
\usepackage{amssymb,amsfonts}
\usepackage{graphicx}                 

\newcommand{\bE}{\mathbf{E}}

\newcommand{\bH}{\mathbf{H}}

\newcommand{\bx}{\mathbf{x}}
\newcommand{\by}{\mathbf{y}}

\newcommand{\Bu}{\boldsymbol{u}}
\newcommand{\Bv}{\boldsymbol{v}}

\newcommand{\cE}{\mathcal{E}}
\def\HH{\mathbf{H}}

\begin{document}


\title{Isoptic curves of conic sections in constant curvature geometries}

\author{G\'eza~Csima and Jen\H o~Szirmai, \\ 
Budapest University of Technology and Economics, \\
Institute of Mathematics,
Department of Geometry Budapest, \\
P.O. Box 91, H-1521 \\
csgeza@math.bme.hu, szirmai@math.bme.hu}
\maketitle
\thanks{}



\begin{abstract}
In this paper we consider the isoptic curves on the 2-dimensional geometries of constant curvature $\bE^2,~\bH^2,~\cE^2$. 
The topic is widely investigated in the Euclidean plane $\bE^2$ see for example \cite{CMM91} and \cite{Wi} and the references given there, but
in the hyperbolic and elliptic plane there are few results in this topic (see \cite{CsSz1} and \cite{CsSz2}). In this paper we 
give a review on the preliminary results of the isoptics of Euclidean and hyperbolic curves and develop a procedure to study the isoptic curves in the 
hyperbolic and elliptic plane geometries and apply it for some geometric objects e.g. proper conic sections. We use for the computations the classical 
models which are based on the projectiv interpretation of the hyperbolic and elliptic geometry and in this manner the isoptic curves 
can be visualized on the Euclidean screen of computer.
\end{abstract}

\newtheorem{theorem}{Theorem}[section]
\newtheorem{corollary}[theorem]{Corollary}
\newtheorem{lemma}[theorem]{Lemma}
\newtheorem{exmple}[theorem]{Example}
\newtheorem{defn}[theorem]{Definition}
\newtheorem{rmrk}[theorem]{Remark}
\newenvironment{definition}{\begin{defn}\normalfont}{\end{defn}}
\newenvironment{remark}{\begin{rmrk}\normalfont}{\end{rmrk}}
\newenvironment{example}{\begin{exmple}\normalfont}{\end{exmple}}
\newtheorem{remarque}{Remark}



\section{Introduction}

Let $X$ be one of the geometries of contant curvature $\bE^2,~\bH^2,~\cE^2$. 
The isoptic curve of an arbitrary given plane curve $\mathcal{C}$ is the locus of points $P \in X$ where $\mathcal{C}$ is seen under
a given fixed angle $\alpha$ $(0<\alpha <\pi)$. 
An isoptic curve formed from the locus of tangents meeting at right angles are called orthoptic curve. 
The name isoptic curve was suggested by Taylor in his work \cite{T} in 1884.

First we consider the Euclidean plane geometry($X=\bE^2$). The easiest case if $\mathcal{C}$ is a line segment then the set of all points (locus) 
for which a line segment can be seen at angle 
$\alpha$ contains two arcs in both halfplane of the line segment, each is with central angle $2\alpha$. 
In the special case of $\alpha=\frac{\pi}{2}$, we get exactly one circle, called Thales circle (without the endpoints of the given segment) 
with center the middle of the line segment. 

In \cite{CMM91} and \cite{CMM96} the isoptic curves of the closed, stricktly convex curves are studied, using their support function. 
The papers \cite{Wu71-1} and \cite{Wu71-2} deal with curves having a
circle or an ellipse for an isoptic curve. Further curves appearing as isoptic curves are well studied in the Euclidean plane geometry 
$\bE^2$, see e.g. \cite{Lo},\cite{Wi}. 
Isoptic curves of conic sections have been studied in \cite{H} and \cite{S}. 
A lot of papers concentrate on the properties of the isoptics e.g. \cite{MM} and \cite{M} and the reference given there. 

There are a lot of possibility, to give the equations of the isoptics of conic sections (see e.g \cite{Lo}), for instance, they can 
be determined by the {\it constuction method of the tangent lines from outher point}. 
We have illustrated this procedure (see \cite{St}) on the following pictures:
\begin{figure}[ht]
\centering
\includegraphics[scale=0.245]{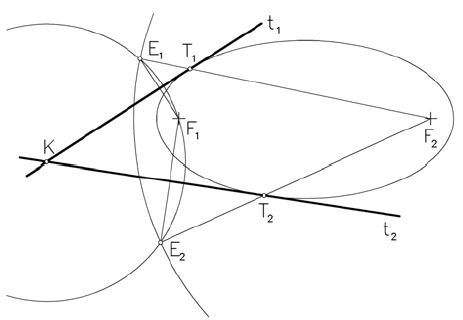}\includegraphics[scale=0.345]{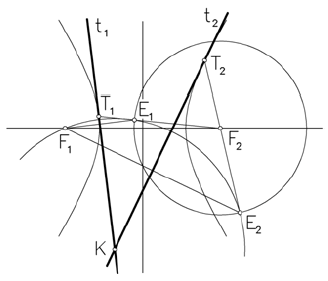}\includegraphics[scale=0.375]{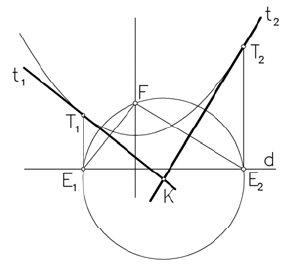}
\caption{Tangent lines from outher point $K$}
\label{fig:erint}
\end{figure}

To get the isoptics, we have to solve equation systems genereted by two circle equations
(ellipse, hyperbola) or a circle and a line equation (parabola) and using the scalar product, 
we have to fix the angle of the tangent lines. 
In case of the hyperbola there is no proper affect point, if the outher point is contained by one of its 
asymptotas, but the asymptotas can be considered as tangent lines according to projective approach. 
From this methode, we get the following equations for the isoptic curves (see \cite{Lo}, as well):
\[\mathnormal{Ellipse:}\ \cos{\alpha}=-\frac{a^2+b^2-x^2-y^2}{\sqrt{\left(-a^2+b^2+x^2\right)^2+2 y^2 \left(a^2-b^2+x^2\right)+y^4}}\]
where the ellipse is given by its equation $\frac{x^2}{a^2}+\frac{y^2}{b^2}=1$,
\[\mathnormal{Hyperbola:}\ \cos^2{\alpha}=\frac{(-a^2+b^2+x^2+y^2)^2}{\left(a^2+b^2-x^2\right)^2+2 y^2 \left(a^2+b^2+x^2\right)+y^4}\]
where the hyperbola is given by its equation $\frac{x^2}{a^2}-\frac{y^2}{b^2}=1$,
\[\mathnormal{Parabola:}\ \cos{\alpha}=-\frac{y}{\sqrt{(p-y)^2+x^2}}\]
where the axis $x$ is the directix, and the focus is in $(0,p)$.
\begin{remark}
\begin{enumerate}
\item
In case of hyperbola, the two asymptotas split the space into four domains, two of them contains a hyperbola branch (focal domains), 
the other ones are empty. Let $P$ an outher point of the hyperbola. 
If $P$ is in a focal domain, then the tangent lines affects the same branch of the hypebola, else they affect both of the branches. 
In these cases, the isoptic angle are complementary to each other, i.e. the sum of them are $\pi$. 
Therfore, we take the squere of the equation thus we obtain both types of isoptic curves.
\item
The numerator is larger then zero for every $(x,y)$, if $b>a$, therefore the isoptic curve do not exists in the interval 
$$\left( \arccos\left(\frac{b^2-a^2}{b^2+a^2}\right),\arccos\left(\frac{a^2-b^2}{a^2+b^2}\right) \right)$$ if the condition $b>a$ holds. Otherwise, the isoptic curve exists for every $\alpha \in (0,\pi)$.
\end{enumerate}
\end{remark}

We have illustrated the isoptic curves to some conic sections in Euclidean plane $\bE^2$ in Figures 2-3:

\begin{figure}[ht]
\centering
\includegraphics[scale=0.61]{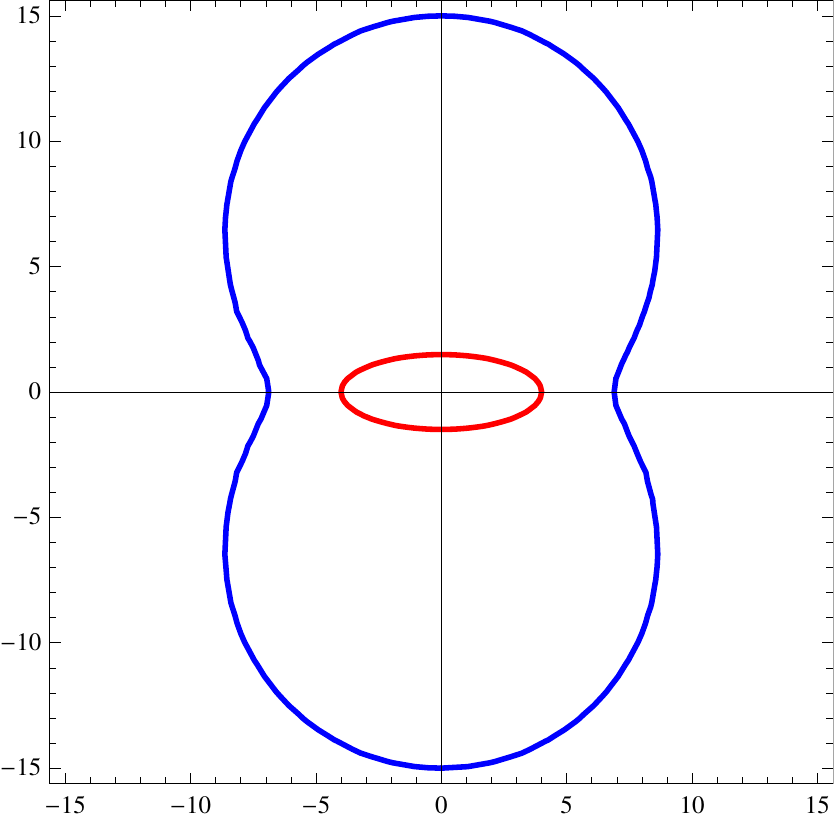} \includegraphics[scale=0.61]{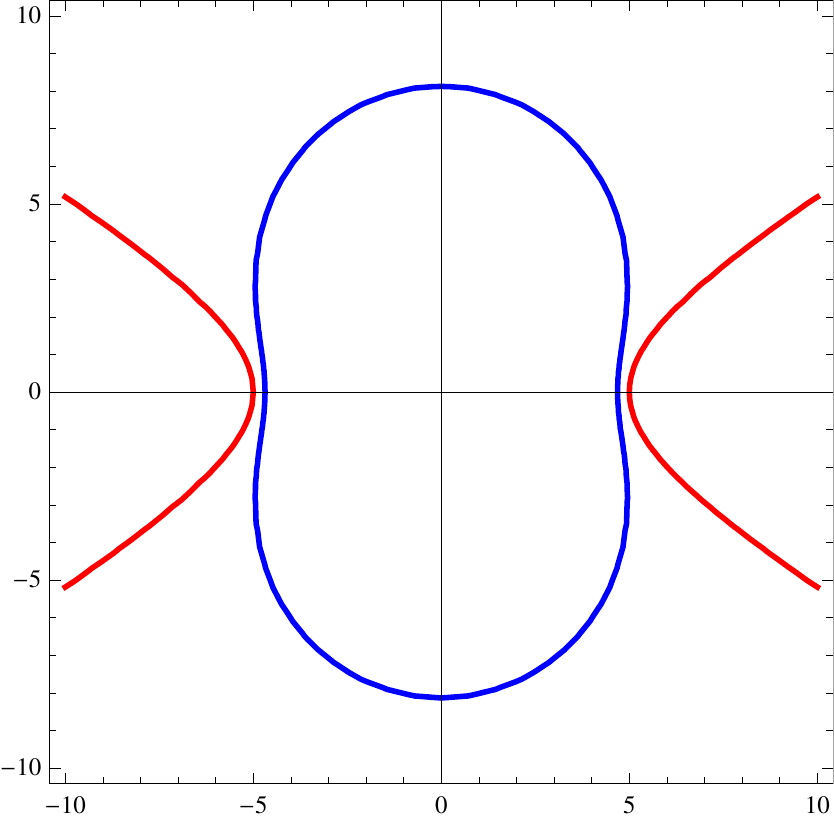}
\caption{Isoptic curve to the Euklidean ellipse (left) and hyperbola (right) with parameters: $a=4$,$b=1,5$, $\alpha=\pi/6$, and $a=5$, $b=3$, $\alpha=\pi/3$}
\label{fig:eukel1}
\end{figure}\begin{figure}[ht]
\centering
\includegraphics[scale=0.61]{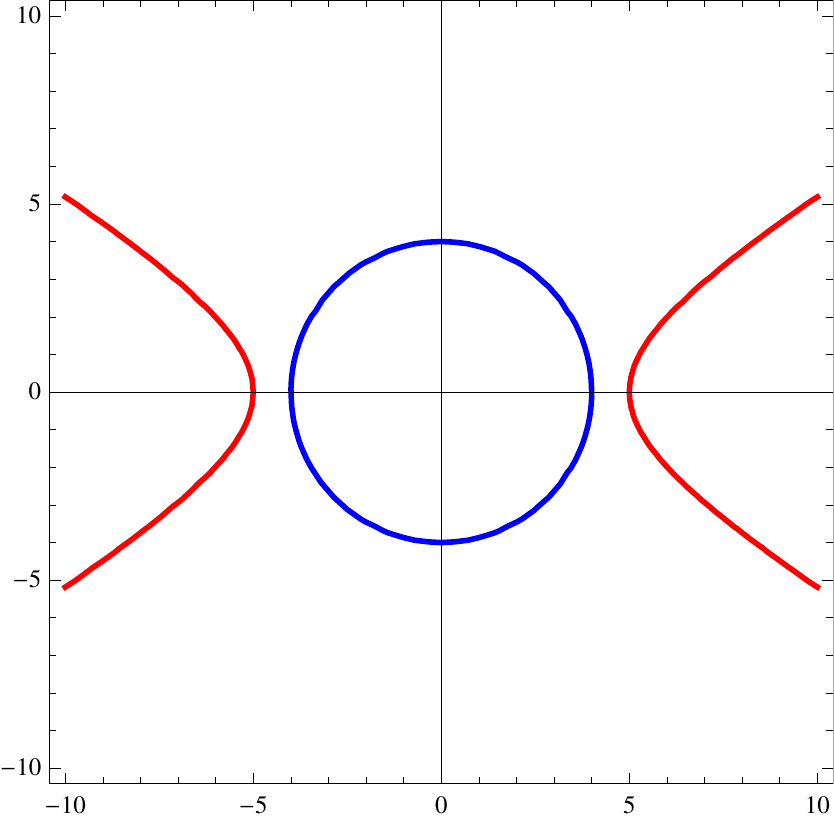} \includegraphics[scale=0.6]{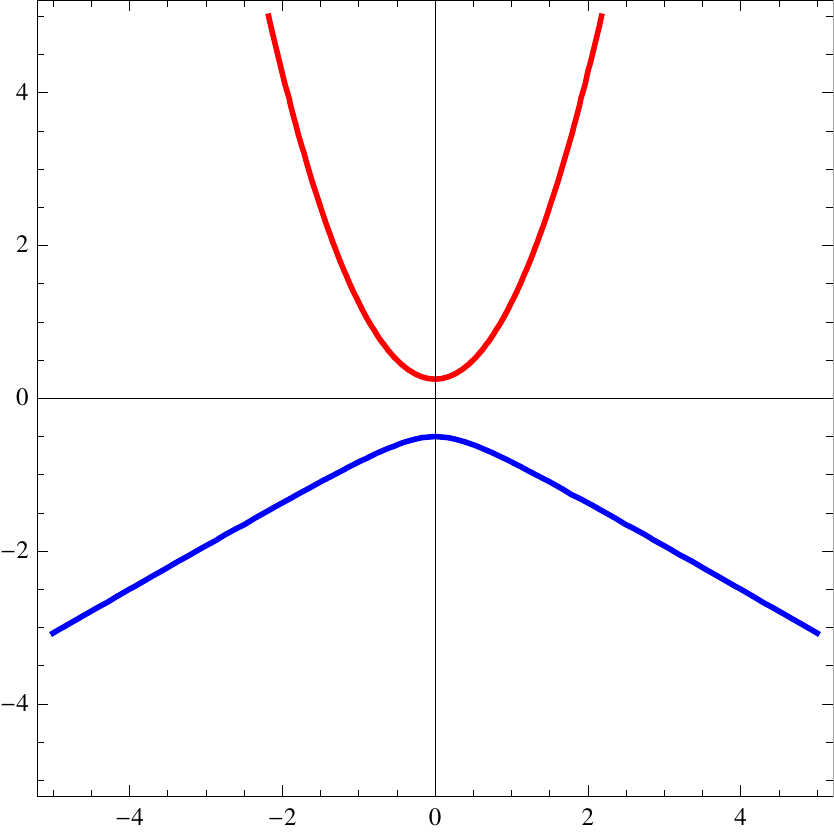}
\begin{center}
\caption{Isoptic curve to the Euklidean hyperbola (left) and parabola (right) with parameters: $a=5$, $b=3$, $\alpha=\pi/2$, and $p=1/2$, $\alpha=\pi/3$}
\end{center}
\label{fig:eukhip2}
\end{figure}
If $X$ is the hyperbolic plane geometry, we have only few results. The isoptic curves of the hyperbolic line segment, ellipses and parabolas are 
determined in \cite{CsSz1} and \cite{CsSz2}. 

As far as we know, there are no results in the elliptic geometry $\cE^2$. 

In this paper, we develop a method, based on the projective interpretation of the hyperbolic and elliptic geometry, 
to determine the isoptic curve of a given plane curve $\mathcal{C}$ and we apply our procedure to the hyperbolic hyperbola with proper foci, 
elliptic line segments and elliptic conic sections, moreover, we visualize them for some angles.  
\section{The projective model}
\subsection{Hyperbolic geometry $\bH^d$}
For the $d$-dimensional hyperbolic space $(d \ge 2)$
$\bH^d$ we use the projective model in Lorentz space
$\bE^{1,d}$ of signature $(1,d)$, i.e.~$\bE^{1,d}$ is
the real vector space $\mathbf{V}^{d+1}$ equipped with the bilinear
form of signature $(1,d)$
\begin{equation}
\bH^d:\ \ \langle ~ \mathbf{x}, ~ \mathbf{y} \rangle = -x^0y^0+x^1y^1+ \dots + x^d y^d \tag{2.1 a}
\end{equation}
where the non-zero vectors
$$
\mathbf{x}=(x^0,x^1,\dots,x^d)\in\mathbf{V}^{d+1} \ \  \text{and} \ \ \mathbf{y}=(y^0,y^1,\dots,y^d)\in\mathbf{V}^{d+1},
$$
are determined up to real factors and they represent points of $\bH^d$ in
$\mathcal{P}^d(\mathbb{R})$. The proper points of $\bH^d$ are represented as the
interior of the absolute quadratic form
\begin{equation}
Q=\{[\mathbf{x}]\in\mathcal{P}^d | \langle ~ \mathbf{x},~\mathbf{x} \rangle =0 \}=\partial \bH^d \tag{2.2}
\end{equation}
in real projective space $\mathcal{P}^d(\mathbf{V}^{d+1},
\mbox{\boldmath$V$}\!_{d+1})$. All proper interior point $[\mathbf{x}] \in \bH^d$ are characterized by
$\langle ~ \mathbf{x},~\mathbf{x} \rangle < 0$.

The points on the boundary $\partial \bH^d $ in
$\mathcal{P}^d$ represent the absolute points at infinity of $\bH^d $.
Points $[\mathbf{y}]$ with $\langle ~ \mathbf{y},~\mathbf{y} \rangle >
0$ lie outside $\partial \bH^d $ and are called outer points
of $\bH^d $. 

Let $P([\mathbf{x}]) \in \mathcal{P}^d$; a point
$[\mathbf{y}] \in \mathcal{P}^d$ is said to be conjugate to
$[\mathbf{x}]$ relative to $Q$ when $\langle ~
\mathbf{x},~\mathbf{y} \rangle =0$. The set of all points conjugate
to $P([\mathbf{x}])$ form a projective (polar) hyperplane
\begin{equation}
pol(P):=\{[\mathbf{y}]\in\mathcal{P}^d | \langle ~ \mathbf{x},~\mathbf{y} \rangle =0 \}. \tag{2.3}
\end{equation}

Hence the bilinear form of $Q$ by (2.1 a) induces a bijection
(linear polarity $\mathbf{V}^{d+1} \rightarrow
\mbox{\boldmath$V$}\!_{d+1})$)
from the points of $\mathcal{P}^d$
onto its hyperplanes.

Point $X [\mathbf{x}]$ and the hyperplane $\alpha
[\mbox{\boldmath$a$}]$ are called incident if the value of the
linear form $\mbox{\boldmath$a$}$ on the vector $\mathbf{x}$ is equal
to zero; i.e., $\mathbf{x}\mbox{\boldmath$a$}=0$ ($\mathbf{x} \in \
\mathbf{V}^{d+1} \setminus \{\mathbf{0}\}, \ \mbox{\boldmath$a$} \in
\mbox{\boldmath$V$}_{d+1} \setminus \{\mbox{$\boldsymbol{0}$}\}$).
Straight lines in $\mathcal{P}^d$ are characterized by the
2-subspaces of $\mathbf{V}^{d+1} \ \text{or $(d-1)$-spaces of} \
\mbox{\boldmath$V$}\!_{d+1}$ .

In this paper we set the sectional curvature of $\bH^d$,
$K=-k^2$, to be $k=1$. The distance $s$ of two proper points
$[\mathbf{x}]$ and $[\mathbf{y}]$ is calculated by the formula:
\begin{equation}
\cosh{{s}}=\frac{-\langle ~ \mathbf{x},~\mathbf{y} \rangle }{\sqrt{\langle ~ \mathbf{x},~\mathbf{x} \rangle
\langle ~ \mathbf{y},~\mathbf{y} \rangle }} . \tag{2.4 a}
\end{equation}

%

Suppose that $[\Bu]$ and $[\Bv]$ are proper lines of $\bH^d$. In $\bH^d$ they intersect in a proper point 
if $\langle \Bu,\Bu \rangle \langle \Bv,\Bv \rangle - \langle \Bu,\Bv \rangle^2 >0$. Their dihedral angle 
$\alpha (\Bu,\Bv)$ can be measured by
\begin{equation}
\cos{\alpha}=\frac{-\langle \Bu,\Bv \rangle}{\sqrt{
\langle \Bu,\Bu\rangle \langle \Bv,\Bv \rangle}}.
\tag{2.5 a}
\end{equation}

\subsection{Elliptic geometry $\cE^d$}

Similarly to the hyperbolic case we can introduce the projective model of the elliptic geometry $\cE^d$ $(d \ge 2)$. 
Let $\bx,\by\in\mathbf{V}^{d+1}$ be given, 
then
\begin{equation}
\cE^d:\ \ \langle ~ \mathbf{x}, ~ \mathbf{y} \rangle = x^0y^0+x^1y^1+x^2y^2 \dots + x^d y^d.  \tag{2.1 b}
\end{equation}
The distance between $\bx$, $\by$ can be measured by the following formula:
\begin{equation}
\cos{{s}}=\frac{\langle ~ \mathbf{x},~\mathbf{y} \rangle }{\sqrt{\langle ~ \mathbf{x},~\mathbf{x} \rangle
\langle ~ \mathbf{y},~\mathbf{y} \rangle }} . \tag{2.4 b}
\end{equation}
And if $\Bu$ and $\Bv$ are lines of the $\cE^d$, then their angle $\alpha$ can be determined:
\begin{equation}
\cos{\alpha}=\frac{\langle \Bu,\Bv \rangle}{\sqrt{
\langle \Bu,\Bu\rangle \langle \Bv,\Bv \rangle}}.
\tag{2.5 b}
\end{equation}

\section{Isoptic curve of the line segment on the hyperbolic and elliptic plane}

In this section, we examine the hyperbolic and elliptic cases together. Let $\epsilon=-1$ if $X$ is the hyperbolic geometry and $\epsilon=1$ if $X$ is the elliptic geometry.

Let two points $A$ and $B$ be given in the plane. We can assume without loss of generality that 
their homogeneous coordinates are $A[{\bf{a}}]\sim(1,a,0)$ and $B[{\bf{b}}]\sim(1,-a,0)$, where $(a\in ]0,1])$. 
We consider two straight lines $u[\Bu]$ and $v[\Bv]$ where $\Bu \sim (1,u_1,u_2)^{T}$ passes trough points $A$ and $P$, and $\Bv \sim (1,v_1,v_2)^{T}$ 
passes through points $B$ and $P$. By the incidence formula we get the following equations:  
\begin{equation}
\begin{gathered}
A \in u \Leftrightarrow (1,a,0)\begin{pmatrix}
1\\
u_1\\
u_2 \\
\end{pmatrix}=0 \ \Leftrightarrow u_1=-\frac{1}{a} \ \\
B \in v \Leftrightarrow (1,-a,0)\begin{pmatrix}
1\\
v_1\\
v_2 \\
\end{pmatrix}=0, \ \ \Leftrightarrow v_1=\frac{1}{a} \tag{3.1}
\end{gathered}
\end{equation}
\begin{equation}
\begin{gathered}
P \in u \leftrightarrow (1,x,y)\begin{pmatrix}
1\\
u_1\\
u_2 \\
\end{pmatrix}=0 \ \Leftrightarrow u_2=-\frac{a-x}{ya}, \ \ y\ne 0, \\
P \in v \Leftrightarrow (1,x,y)\begin{pmatrix}
1\\
v_1\\
v_2 \\
\end{pmatrix}=0, \ \ \Leftrightarrow v_2=-\frac{a+x}{ya}, \ \ y\ne0. \tag{3.2}
\end{gathered}
\end{equation}

The angle $\alpha$ between the above straigt lines can be determined by the formula (2.5 a) and (2.5 b):
\[
\cos(\alpha)=\frac{\epsilon(\epsilon+u_1v_1+u_2v_2)}{\sqrt{(\epsilon+u_1^2+u_2^2)(\epsilon+v_1^2+v_2^2)}}. \notag
\]
Substituting coordinates from (3.1) and (3.2) into the above equation, we obtain:
\begin{theorem} 
Let suppose that a line segment is given by $A[{\bf{a}}]\sim(1,a,0)$ and $B[{\bf{b}}]\sim(1,-a,0)$. 
Then for a given $\alpha(0<\alpha<\pi)$, the $\alpha$-isoptic curve of the $AB$ on the hyperbolic and elliptic plane has
an equation of the form:
\begin{equation}
\cos(\alpha)=\frac{\epsilon(\epsilon-\frac{1}{a^2}+\frac{a^2-x^2}{y^2a^2})}{\sqrt{(\epsilon+\frac{1}{a^2}+(\frac{a-x}{ya})^2)(\epsilon+\frac{1}{a^2}+(\frac{a+x}{ya})^2)}},
\label{hipszog} \tag{3.3}
\end{equation}
where $\epsilon=-1$ in the hyperbolic geometry and $\epsilon=1$ in the elliptic geometry.
\end{theorem}

\begin{remark}
\begin{enumerate}
\item In the hyperbolic case we obtain the orthoptic curve $\mathcal{G}_{\pi/2}$ if $\alpha=\pi/2$ with equation:
\begin{equation}
\frac{x^2}{a^2}+\frac{y^2}{\frac{a^2}{1+a^2}}=1. \tag{3.4}
\end{equation}
This is an ellipse (without endpoints of the given segment) in the Euclidean sense, and it can be called the 
{\it Thales curve} in the hyperbolic geometry.

If we increase the parameter $a$ then the {\it Thales curve} tends to a hypercycle (or equidistant curve) in the hyperbolic plane, that means
the hypercycle is a special type of the orthoptic curves:
$$
x^2+2y^2=1. 
$$
\item In the hyperbolic plane, if $a$ tends to $1$, than the (3.3) equation tends to a the following equation:
\begin{equation}
x^2+\left(\frac{y}{\cos\left(\frac{\alpha}{2}\right)}\right)^2=1.
\notag
\end{equation}
\end{enumerate}
\end{remark}

\begin{figure}[ht]
\centering
\includegraphics[scale=0.615]{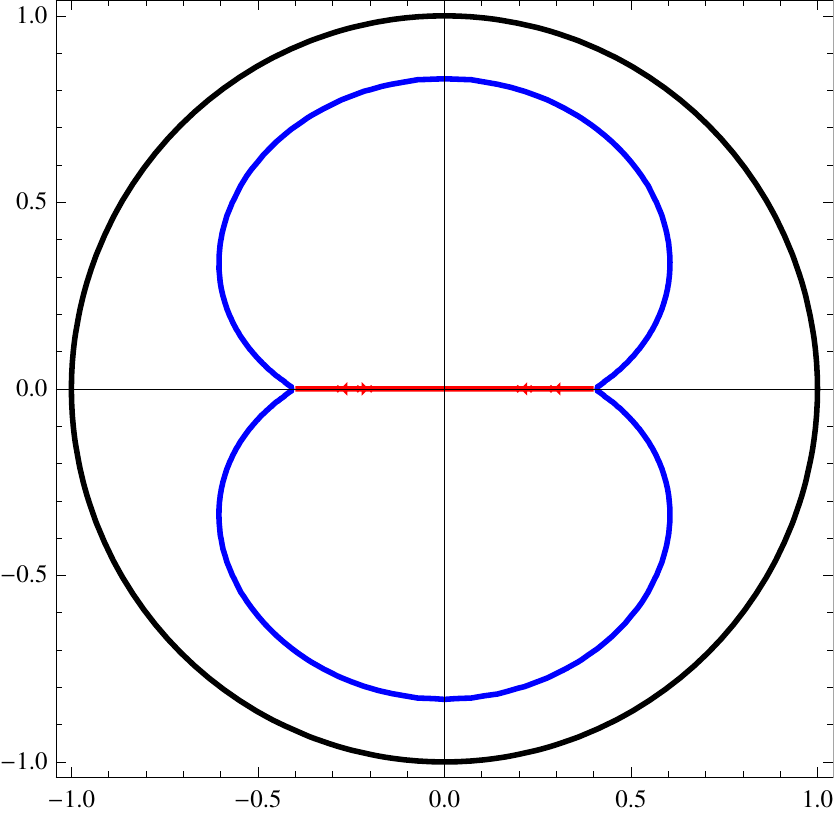} \includegraphics[scale=0.615]{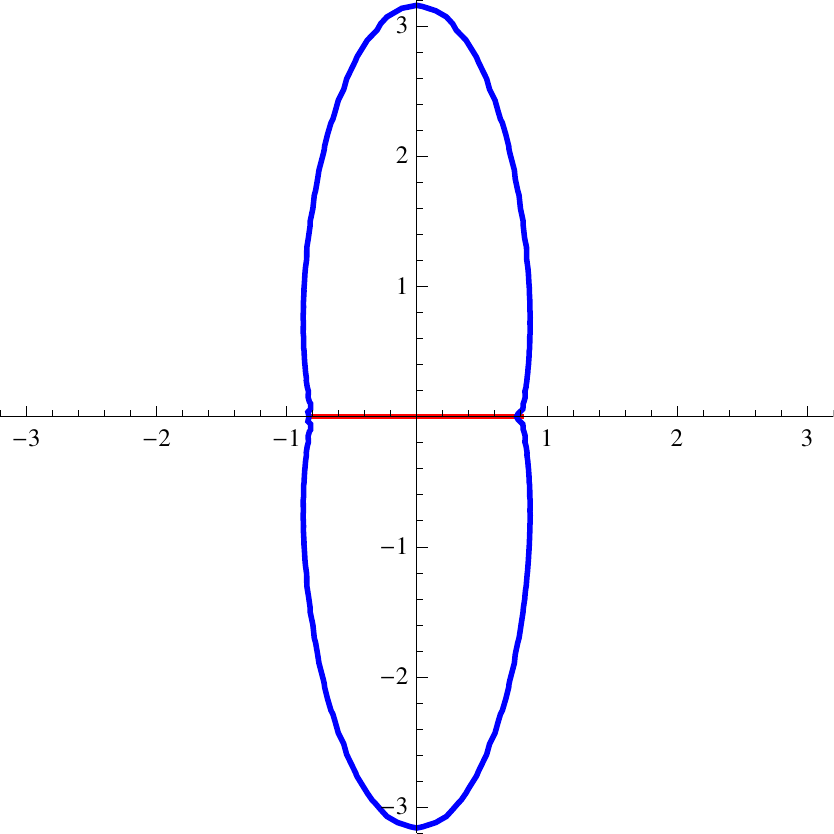}
\caption{The isoptic curve to hyperbolic (left) and elliptic (rigth) line segment with
parameters: Hyperbolic: $a=0.4$, $\alpha=\pi/6$,~ ~Elliptic: $a=0.8$,  $\alpha=8\pi/18$}
\end{figure}

\section{The general method}
The method bellow has been implicated by the isoptic curve of the line segment. That procedure can be used to develop a more general 
method to determine the isoptic curves. 
Let a conic section $C$ and one of its point $P$ be given.
Using the implicit function theory and the equation of $C$, 
we can determine the equation of the tangent line in this point. 
After that, there can be given an equation system to the coordinates of the tangent point from an external point $K$. 
This point have to satisfy the equation of the given curve and the tangent lines to this points have to contain $K$. 
This system can be solved for every $K(1,x_0,y_0)$ outher point with respect to the parameters $x_0$, $y_0$. 
In general case, the solutions (both of them) are complicated. 
Now, we have to follow the upper method using the coordinates of the tangent points. The equation of the tangent lines from $K$ 
can be determined by solving an equation system for its coordinates.
Finally, we have to fix the angle of the straight lines, and we get the equation of those point, from the given curve can be seen under the 
given $\alpha$ $(0<\alpha<\pi)$, which is the isoptic curve.
\subsection{On the hyperbolic plane}
\subsubsection{Equation of the hyperbolic ellipse and hyperbola}
Now, we define the {\it proper central conic sections} and give their equations. 
\begin{defn}
The proper hyperbolic ellipse is the locus of all points of the hyperbolic plane whose distances to two proper fixed points add 
to the same constant ($2a$).
\end{defn}
\begin{defn}
The proper hyperbolic hyperbola is the locus of points where the absolute value of the difference of the distances to the two proper foci 
is a constant ($2a$).
\end{defn}

We discuss the ellipse and the hyperbola together. We can suppose that the two foci are equidistant from the orign $O$, 
both fits on the axis $x$ with coordinates $F_1[\textbf{f}_{1}]\sim(1,f,0)$ and $F_2[\textbf{f}_{2}]\sim(1,-f,0)$ where $0<f<1$. 
Let $P[\textbf{p}] \sim (1,x,y) \in \HH^2$ a point of the conic section. Using (2.4 a) we obtain the following equation:
\begin{equation}
\epsilon_1\cosh ^{-1}\left(\frac{-\left\langle \textbf{p},\textbf{f}_{1}\right\rangle}{\sqrt{\left\langle \textbf{p},\textbf{p}\right\rangle \left\langle 
\textbf{f}_{1},\textbf{f}_{1}\right\rangle}}\right)+\epsilon_2\cosh ^{-1}\left(\frac{-\left\langle \textbf{p},\textbf{f}_{2}\right\rangle}{\sqrt{\left\langle 
\textbf{p},\textbf{p}\right\rangle \left\langle \textbf{f}_{2},\textbf{f}_{2}\right\rangle}}\right)=2a \Leftrightarrow \notag
\end{equation}
\begin{equation*}
\begin{gathered}
\Leftrightarrow \epsilon_2\cosh ^{-1}\left(\frac{-(-1-xf)}{\sqrt{(-1+x^2+y^2)(-1+f^2)}}\right)=\\
= 2a-\epsilon_1\cosh ^{-1}\left(\frac{-(-1+xf)}{\sqrt{(-1+x^2+y^2)(-1+f^2)}}\right),
\nonumber
\end{gathered}
\end{equation*}
where $\epsilon_{1,2}=\pm1$ and $\epsilon_1+\epsilon_2\geq0$.
Taking the $\cosh(\ )$ function for both sides, we get the following equations:
\begin{equation}
\begin{gathered}
\frac{1+xf}{\sqrt{(-1+x^2+y^2)(-1+f^2)}}=
\cosh(2a)\frac{1-xf}{\sqrt{(-1+x^2+y^2)(-1+f^2)}}-\\-\epsilon_1
\sinh(2a)\sinh\left(\cosh^{-1}\left(\frac{1-xf}{\sqrt{(-1+x^2+y^2)(-1+f^2)}}\right)\right). \notag
\end{gathered}
\end{equation}
The next equation is obtained by applying the formula $$\sinh\left(\cosh^{-1}(t)\right)=\sqrt{\cosh^{2}\left(\cosh^{-1}(t) \right)-1}=\sqrt{t^{2}-1}$$ and
by multiplying both sides by $\sqrt{(-1+x^2+y^2)(-1+f^2)}$:  

\begin{equation}
\begin{gathered}
1+xf=\cosh(2a)(1-xf)-\epsilon_1\sinh(2a)\sqrt{(1-xf)^2-(1-x^2-y^2)(1-f^2)}.
\notag
\end{gathered}
\end{equation}

Now, if we simplify this equation, and take its squere, there is no $\epsilon$ in it. Finally, we get the following equation:

\begin{equation}
\left(\frac{x}{\tanh(a)}\right)^2+\frac{y^2}{1+\frac{1}{(f^2-1)\cosh^2(a)}}=1. \tag{5.1}
\end{equation}

If the distance between the two foci is lesser than $2a$, it is an ellipse, if larger than it is a hyperbola.

\begin{equation}
\begin{gathered}
2a<>d(F_1,F_2)\Leftrightarrow\cosh(2a)<>\cosh(d(F_1,F_2))=\frac{-\left\langle F_1,F_2\right\rangle}{\sqrt{\left\langle F_1,F_1\right\rangle\left\langle F_2,F_2\right\rangle}}=\\
=\frac{1+f^2}{1-f^2}\Leftrightarrow 2\cosh^2(a)-1<>\frac{2}{1-f^2}-1\Leftrightarrow1<>\frac{1}{\cosh^2(a)(1-f^2)}\Leftrightarrow \\
\Leftrightarrow 1+\frac{1}{\cosh^2(a)(f^2-1)}<>0
\notag
\end{gathered}
\end{equation}

Therefore, the hyperbolic ellipse and hyperbola are also ellipse and hyperbola in the model.
\subsubsection{Isoptic curve of the hyperbolic ellipse and hyperbola}
Now, we will use the about discribed method to determine the isoptic curves to hyperbolic ellipses and hyperbolas. 

The firs step is to determine the equation of the tangent lines ($y\neq0$):
\begin{equation}
y'=-\frac{x}{y}\left(1+\frac{f^{2}}{\sinh^{2}(a)(f^2-1)}\right).
\tag{5.2}
\end{equation}
The equation above goes continously to $x^2=\tanh^2(a)$, if $y\rightarrow0$.

After that, we have to solve the following equation system for $\tilde{X}_i[\mathbf{\tilde{x}_i}]\sim(1,\tilde{x}_i,\tilde{y}_i)$, $(i=1,2)$, where $X[\mathbf{x}]\sim(1,x,y)$ is a point in the Cayley-Klein model:
\begin{equation}
\begin{gathered}
y=-\frac{\tilde{x}}{\tilde{y}}\left(1+\frac{f^{2}}{\sinh^{2}(a)(f^2-1)}\right)(x-\tilde{x})+\tilde{y} \\
\left(\frac{\tilde{x}}{\tanh(a)}\right)^2+\frac{\tilde{y}^2}{1+\frac{1}{(f^2-1)\cosh^2(a)}}=1.
\tag{5.3}
\end{gathered}
\end{equation}

It is not so hard, to determine the roots, but because of the complexity of the result, we ignore it. 
Now we need $\Bu \sim (1,u_1,u_2)^T$ and 
$\Bv \sim (1,v_1,v_2)^T$ straigth lines, fits on respectively $P$, $\tilde{X_1}$ and $P$, $\tilde{X_2}$. 

Using that $(\tilde{x}_{i},\tilde{y}_{i})$ is known, we get the following equation system:
\begin{equation}
\begin{gathered}
1+u_{1}x+u_{2}y=0\\
1+u_{1}\tilde{x}_1+u_{2}\tilde{y}_1=0,
\tag{5.4}
\end{gathered}
\end{equation}
\begin{equation}
\begin{gathered}
1+v_{1}x+v_{2}y=0\\
1+v_{1}\tilde{x}_2+v_{2}\tilde{y}_2=0.
\tag{5.5}
\end{gathered}
\end{equation}

Solving these systems, we can determine the tangent lines, for all $X[\textbf{x}]$ which is outside of the conic section. The last step in the method of the previous section is to fix the angle. We summarize our results in the following theorem:

\begin{theorem}
Let a hyperbolic ellipse or hyperbola be centered at the origin in the projective model given by its semimajor axis $a$ and foci
$F_1[\textbf{f}_{1}]\sim(1;f;0)$, 
$F_2[\textbf{f}_{2}]\sim(1;-f;0)$, $(0<f<1)$ such that $2a>d(F_1,F_2)$ or $2a<d(F_1,F_2)$ holds. 
The $\alpha$ and $(\pi-\alpha)$-isoptic curves $(0<\alpha<\pi)$ of the considered ellipse or hyperbola in the hyperbolic plane has
the equation:
\begin{equation}
\cos^2 (\alpha )=\frac{\left(\left(f^2-1\right) \cosh (2 a) \left(x^2+y^2-1\right)+f^2 x^2-1\right)^2}{-2 \left(f^2-1\right) y^2\left(f^2+x^2\right)+\left(f^2-x^2\right)^2+\left(f^2-1\right)^2 y^4},
\tag{5.6}
\end{equation}
where $x^2+y^2\leq1$ condition holds.
\end{theorem}

\begin{remark}
\begin{enumerate}
\item The orthoptic curve of the hyperbolic hyperbola exists, if $f\leq\sqrt{1-\frac{1}{\cosh(2a)}}$ and it is an ellipse, similary to the hyperbolic ellipse, with the following equation:
\begin{equation}
\left(1-f^2\right) \cosh (2 a) \left(-1+x^2+y^2\right)+f^2 x^2=1.
\notag
\end{equation}
\item
We note here without any calculations, that the isoptic curve of the hyperbolic ellipse do not exists in the following interval:
$$\alpha\in\left(\arccos\left(\frac{\left(f^2-1\right) \cosh (2 a)+1}{f^2}\right),\arccos\left(\frac{\left(1-f^2\right) \cosh (2 a)-1}{f^2}\right)\right),
$$
if $\left(\frac{1}{1-f^2}+1\right) \text{sech}^2(a)>2$.
\end{enumerate}
\end{remark}

\begin{figure}[ht]
\centering
\includegraphics[scale=0.65]{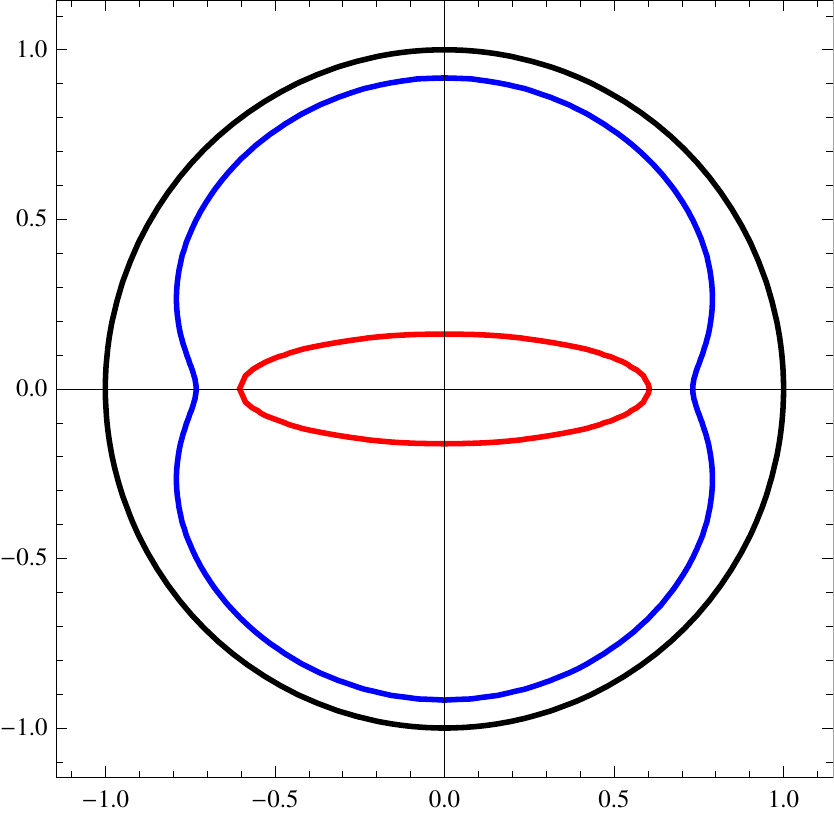} \includegraphics[scale=0.65]{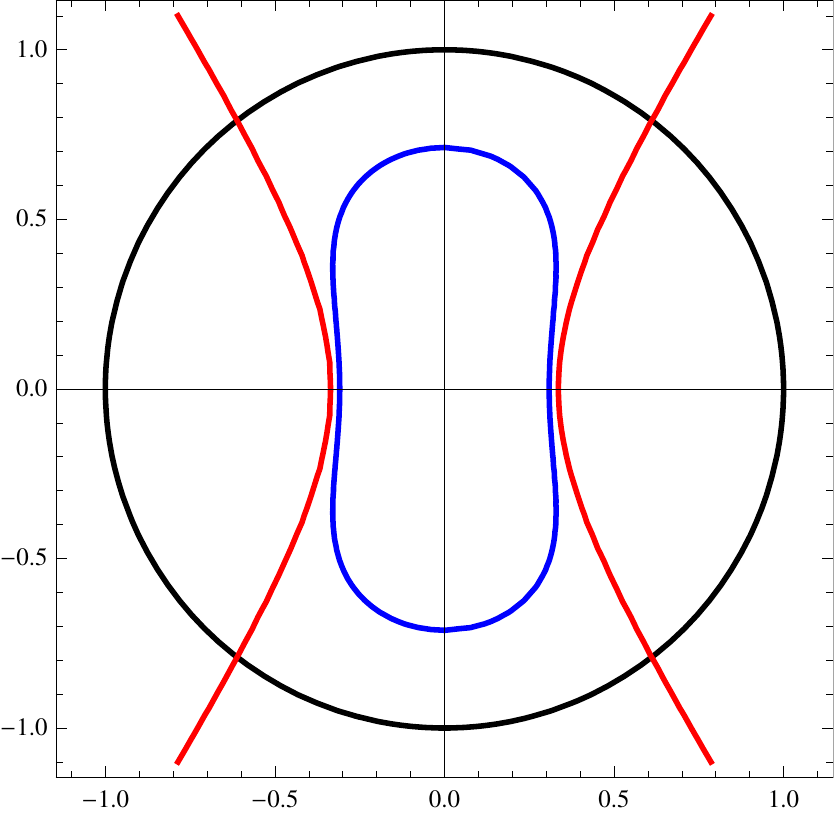}
\caption{The isoptic curve to hyperbolic ellipse (left) and hyperbola (rigth) with parameters: \newline Ellipse: $a=0.7$, $f=0.59$, $\alpha=\pi/6$,~ ~ ~ Hyperbola: $a=0.35$, $f=0.55$, $\alpha=\pi/6$}
\end{figure}

\subsubsection{The equation and the isoptic curve of the hyperbolic parabola}

In this section we define the {\it proper hyperbolic parabolas} and give their equations. 
\begin{defn}
The hyperbolic proper parabola is the set of points ($X[\textbf{x}]\sim(1;x;y) \in\bH^2$) in the hyperbolic plane that are equidistant 
from a proper point (the focus $F$) 
and a proper line (the directrix $e$.) ($s=d(X;F)=d(X;e)$)
\end{defn}

We can assume without loss of generality that the directrix ($e$) is the axis $x$ and the coordinates of the focus point are 
$F[\textbf{f}]\sim(1;0;p)$.

We remark, that the coordinates of the foot point $X'$ of the perpendicular dropped form $X$ into the $x$ axis are $X'[\textbf{x'}]\sim(1;x;0)$. 
Using the (2.4 a) formula:
\begin{equation}
\cosh(s)=\frac{1-py}{\sqrt{1-x^2-y^2}\sqrt{1-p^2}}=\frac{1-x^2}{\sqrt{1-x^2-y^2}\sqrt{1-x^2}}.
\tag{5.7}
\label{hparabeq}
\end{equation}
The equation of the proper hyperbolic
parabola is obtained by (\ref{hparabeq}):%
\begin{equation}
x^2+\frac{(1-py)^2}{1-p^2}=1.
\tag{5.8}
\end{equation}

Using the above methode, we have to solve the folloving equation system for coordinates $\tilde{x},\tilde{y}$:
\begin{equation}
\begin{gathered}
y=\frac{\tilde{x}(1-p^2)}{p-\tilde{y}p^2}(x-\tilde{x})+\tilde{y} \\
x^2+\frac{(1-py)^2}{1-p^2}=1.
\tag{5.9}
\end{gathered}
\end{equation}

In accordance with the method, we have to solve the equation systems (5.4), (5.5). Now, the following theorem hold:
\begin{theorem}
Let a proper hyperbolic parabola be given by its focus $F[\textbf{f}]\sim(1;0;p)$ and its directrix $e$ which is coincide with $x$ axis 
in the projective model. The $\alpha$-isoptic curve of this parabola $(0<\alpha<\pi)$ in the hyperbolic plane have the equation:
\begin{equation}
\cos (\alpha )=\frac{y(py-1)}{\sqrt{(x^2-1)((p^2(x^2-1)+2py+y^2-x^2)}}.
\tag{5.10}
\end{equation}
\end{theorem}

\begin{remark}
The orthoptic curve of the hyperbolic parabola contains two\\ straight line; $y=0$ and $y=\frac{1}{p}$.
\end{remark}

\begin{figure}[ht]
\centering
\includegraphics[scale=0.65]{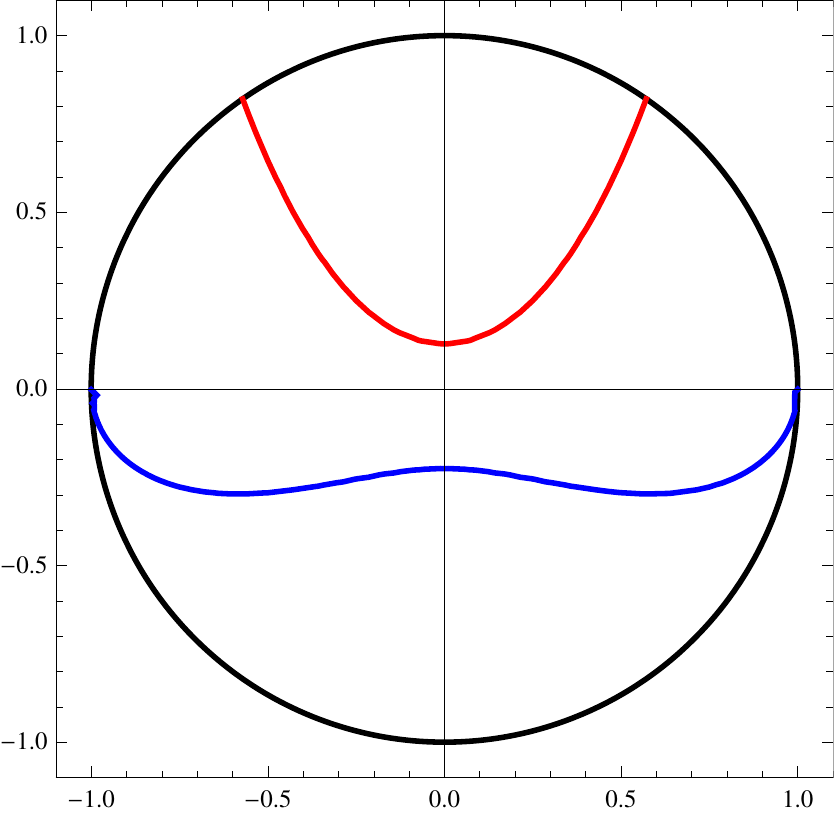} \includegraphics[scale=0.65]{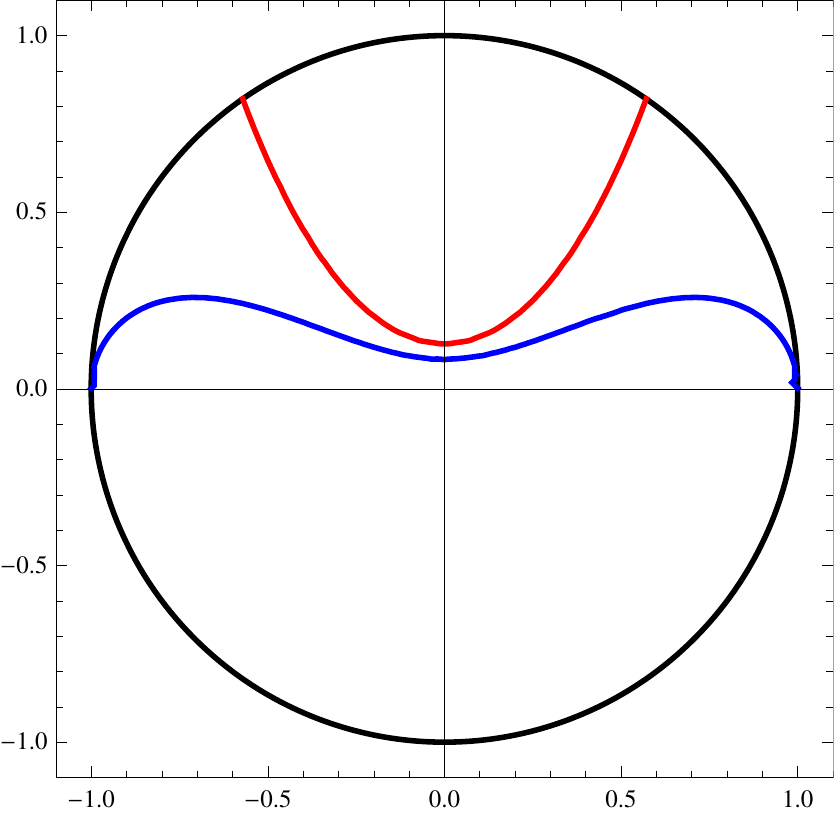}
\caption{The isoptic curve to hyperbolic parabola with parameters: \newline  Parabola: $p=0.25$, $\alpha=\pi/3$,~ ~ ~ ~ ~ ~ ~ ~ ~ ~$p=0.25$, $\alpha=2\pi/3$}
\end{figure}

\subsection{On the elliptic plane}

In this section we will discuss the equations of the conic sections and their isoptics in the elliptic plane $\cE^2$. 

We remark, that on the elliptic plane the maximum distance between two point is less then equal to $\frac{\pi}{2}$, therefore, 
in some cases the given curve cannot be seen under 
arbitrarily small angle.  
\subsubsection{Equation of the elliptic ellipse and hyperbola}

We will follow the deduction process detailed in the previus section to determine the equation of the elliptic ellipse and hyperbola, 
having the following definitions:
\begin{defn}
The elliptic ellipse is the locus of all points of the elliptic plane whose distances to two fixed points add to the same constant ($2a$).
\end{defn}
\begin{defn}
The elliptic hyperbola is the locus of points where the absolute value of the difference of the distances to the two foci is a constant ($2a$).
\end{defn}

Let us suppose that the two foci are equidistant from the orign $O$, both fits on the axis $x$. Then their coordinates are 
$F_1[\textbf{f}_{1}]\sim(1,f,0)$ and $F_2[\textbf{f}_{2}]\sim(1,-f,0)$ where $(0<f<\pi/2)$. 
Let $P[\textbf{p}] \sim (1,x,y) \in \cE^2$ a point of the conic section.

Using (2.4 b) we obtain the following equation:
\begin{equation}
\epsilon_1\cos ^{-1}\left(\frac{\left\langle \textbf{p},\textbf{f}_{1}\right\rangle}{\sqrt{\left\langle \textbf{p},\textbf{p}\right\rangle \left\langle 
\textbf{f}_{1},\textbf{f}_{1}\right\rangle}}\right)+\epsilon_2\cos ^{-1}\left(\frac{\left\langle \textbf{p},\textbf{f}_{2}\right\rangle}{\sqrt{\left\langle 
\textbf{p},\textbf{p}\right\rangle \left\langle \textbf{f}_{2},\textbf{f}_{2}\right\rangle}}\right)=2a \Leftrightarrow \notag
\end{equation}
\begin{equation*}
\begin{gathered}
\Leftrightarrow \epsilon_2\cos ^{-1}\left(\frac{(1-xf)}{\sqrt{(1+x^2+y^2)(1+f^2)}}\right)=\\
= 2a-\epsilon_1\cos ^{-1}\left(\frac{(1+xf)}{\sqrt{(1+x^2+y^2)(1+f^2)}}\right),
\nonumber
\end{gathered}
\end{equation*}
where $\epsilon_{1,2}=\pm1$ and $\epsilon_1+\epsilon_2\geq0$.
Taking the $\cos(\ )$ function for both sides, we get the following equations:
\begin{equation}
\begin{gathered}
\frac{1-xf}{\sqrt{(1+x^2+y^2)(1+f^2)}}=
\cos(2a)\frac{1+xf}{\sqrt{(1+x^2+y^2)(1+f^2)}}+\\+\epsilon_1
\sin(2a)\sin\left(\cos^{-1}\left(\frac{1+xf}{\sqrt{(1+x^2+y^2)(1+f^2)}}\right)\right). \notag
\end{gathered}
\end{equation}
The next equation is obtained by applying the formula $$\sin\left(\cos^{-1}(t)\right)=\sqrt{1-\cos^{2}\left(\cos^{-1}(t) \right)}=\sqrt{1-t^{2}}$$ and
by multiplying both sides by $\sqrt{(1+x^2+y^2)(1+f^2)}$:  

\begin{equation}
\begin{gathered}
1-xf=\cos(2a)(1+xf)+\epsilon_1\sin(2a)\sqrt{(1+x^2+y^2)(1+f^2)-(1+xf)^2}.\notag
\end{gathered}
\end{equation}

Now, if we sort this equation, and take its squere, there is no $\epsilon$ in it. Finally, we get the following equation:

\begin{equation}
\left(\frac{x}{\tan(a)}\right)^2+\frac{y^2}{\frac{1}{(1+f^2)\cos^2(a)}-1}=1. \tag{6.1}
\end{equation}

If the distance between the two foci is lesser than $2a$, it is an ellipse, else it is a hyperbola.

\begin{equation}
\begin{gathered}
2a<>d(F_1,F_2)\Leftrightarrow\cos(2a)<>\cos(d(F_1,F_2))=\frac{\left\langle F_1,F_2\right\rangle}{\sqrt{\left\langle F_1,F_1\right\rangle\left\langle F_2,F_2\right\rangle}}=\\
=\frac{1-f^2}{1+f^2}\Leftrightarrow 2\cos^2(a)-1<>\frac{2}{1+f^2}-1\Leftrightarrow1<>\frac{1}{\cos^2(a)(1+f^2)}\Leftrightarrow \\
\Leftrightarrow 1-\frac{1}{\cos^2(a)(1+f^2)}<>0
\notag
\end{gathered}
\end{equation}

The above first step is a equivalence, because the function $\cos(\ )$ is strictly decreasing in the interval $(0,\pi)$.

We have to take the implicit derivative of (6.1), and solve the equation system for the tangent points 
$\tilde{X}_i[\mathbf{\tilde{x}_i}]\sim(1,\tilde{x}_i,\tilde{y}_i)$, $(i=1,2)$:
\begin{equation}
\begin{gathered}
y=-\frac{\tilde{x}}{\tilde{y}}\left(1-\frac{f^{2}}{\sin^{2}(a)(1+f^2)}\right)(x-\tilde{x})+\tilde{y} \\
\left(\frac{\tilde{x}}{\tan(a)}\right)^2+\frac{\tilde{y}^2}{\frac{1}{(1+f^2)\cos^2(a)}-1}=1.
\tag{6.2}
\end{gathered}
\end{equation}

By solving the equation systems (5.4), (5.5) to these roots $\tilde{X}_{1,2}$ we get the following.
\begin{theorem}
Let an elliptic ellipse or hyperbola be centered at the origin of the projective model, given by its semimajor axis $a$ and its foci 
$F_1[\textbf{f}_{1}]\sim(1;f;0)$, 
$F_2[\textbf{f}_{2}]\sim(1;-f;0)$, $(0<f<1)$ such that $2a>d(F_1,F_2)$ or $2a<d(F_1,F_2)$ holds. 
The $\alpha$ and $(\pi-\alpha)$-isoptic curves $(0<\alpha<\pi)$) of the considered ellipse or hyperbola in the elliptic plane have
the equation:
\begin{equation}
\cos^2 (\alpha )=\frac{\left(\left(1+f^2\right) \cos (2 a) \left(x^2+y^2+1\right)+f^2 x^2-1\right)^2}{2 \left(1+f^2\right) y^2\left(f^2+x^2\right)+\left(f^2-x^2\right)^2+\left(1+f^2\right)^2 y^4}.
\tag{6.3}
\end{equation}
\end{theorem}

\begin{remark}
\begin{enumerate}
\item The orthoptic curve of the elliptic ellipse and hyperbola is an ellipse, with the following equation:
\begin{equation}
\left(1+f^2\right) \cos (2 a) \left(x^2+y^2+1\right)+f^2 x^2=1.
\notag
\end{equation}
\item Ignoring the discussion, the isoptic curve of the elliptic hyperbola exists, if the following formula is True:
\[
\begin{gathered}
\left(\cos\alpha\leq \max\left(\frac{1-(1+f^2)\cos(2a)}{f^2},f^2+(1+f^2)\cos(2a)\right)\right)\wedge\\
\wedge\left(a\geq\frac{\pi}{6}\vee\left(f\leq\sqrt{\frac{1}{\cos(2a)}-1}\right)\vee
\left(\alpha\notin I\right)\right),
\end{gathered} 
\]
where
\[
I=\left(\arccos\left(\frac{(1+f^2)\cos(2a)-1}{f^2}\right),\arccos\left(\frac{1-(1+f^2)\cos(2a)}{f^2}\right)\right).
\]
\end{enumerate}
\end{remark}

\begin{figure}[ht]
\centering
\includegraphics[scale=0.65]{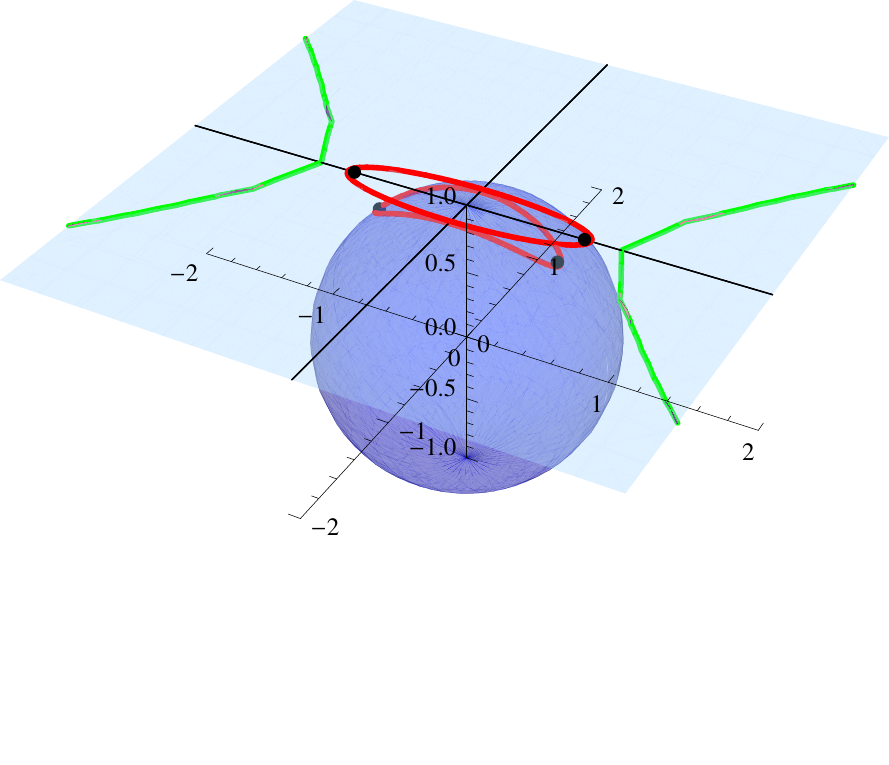} \includegraphics[scale=0.65]{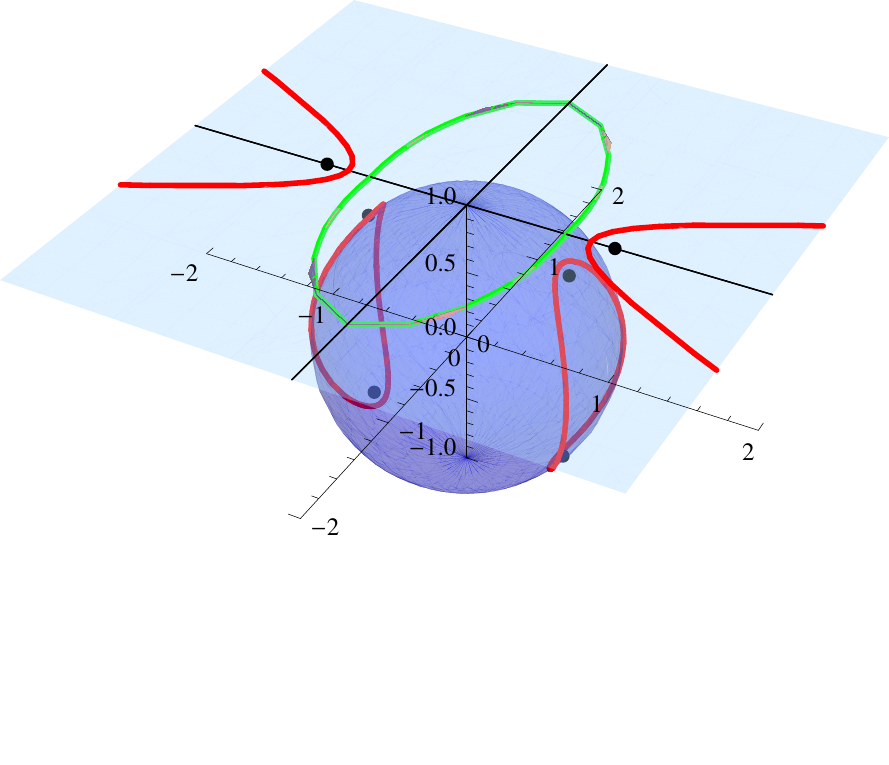}
\caption{The isoptic curve to elliptic ellipse (left) and hyperbola (rigth) with parameters: \newline Ellipse: $a=0.7$, $f=0.8$; $\alpha=\pi/3$~ ~ Hyperbola: $a=0.7$, $f=1$; $\alpha=\pi/2$}
\end{figure}

\subsubsection{The equation and the isoptic curve of the elliptic parabola}

Similarly to the hyperbolic case we define the {\it elliptic parabolas} and give their equations. 
\begin{defn}
A elliptic parabola is the set of points ($X[\textbf{x}]\sim(1;x;y)\in\cE^2$) in the elliptic plane that are equidistant 
from a proper point (the focus $F$) 
and a proper line (the directrix $e$) ($s=d(X;F)=d(X;e)$).
\end{defn}

We can assume without loss of generality that the directrix($e$) is the axis $x$ and the coordinates of the focus point are 
$F[\textbf{f}]\sim(1;0;p)$.
The distances $d(X;F)$ and $d(X;e)$ can be computed by the formulas (2.5 b): 
\begin{equation}
\cos(s)=\frac{1+py}{\sqrt{1+x^2+y^2}\sqrt{1+p^2}}=\frac{1+x^2}{\sqrt{1+x^2+y^2}\sqrt{1+x^2}}.
\tag{6.4}
\label{eparabeq}
\end{equation}
From (\ref{eparabeq}), we obtain the equation of the elliptic parabola:
\begin{equation}
-x^2+\frac{(1+py)^2}{1+p^2}=1.
\tag{6.5}
\end{equation}
Using the same methode, we can solve the foloving equation system for $(\tilde{x},\tilde{y})$:
\begin{equation}
\begin{gathered}
y=\frac{\tilde{x}(1+p^2)}{p+\tilde{y}p^2}(x-\tilde{x})+\tilde{y} \\
-x^2+\frac{(1+py)^2}{1+p^2}=1.
\tag{6.6}
\end{gathered}
\end{equation}

Finally, we have to solve the equation systems (5.4),(5.5). Now, the following theorem hold:
\begin{theorem}
Let a elliptic parabola be given by its focus $F[\textbf{f}]\sim(1;0;p)$ and its directrix $e$ which is coincide with $x$ axis 
in the projective model. The $\alpha$-isoptic curve of this parabola $(0<\alpha<\pi)$ in the elliptic plane has the equation:
\begin{equation}
\cos (\alpha )=\frac{y(py+1)}{\sqrt{(x^2+1)((p^2(x^2+1)-2py+y^2+x^2)}}.
\tag{6.7}
\end{equation}
\end{theorem}

\begin{remark}
The orthoptic curve of the elliptic parabola contains two straight line; $y=0$ and $y=-\frac{1}{p}$.
\end{remark}

\begin{figure}[ht]
\centering
\includegraphics[scale=0.55]{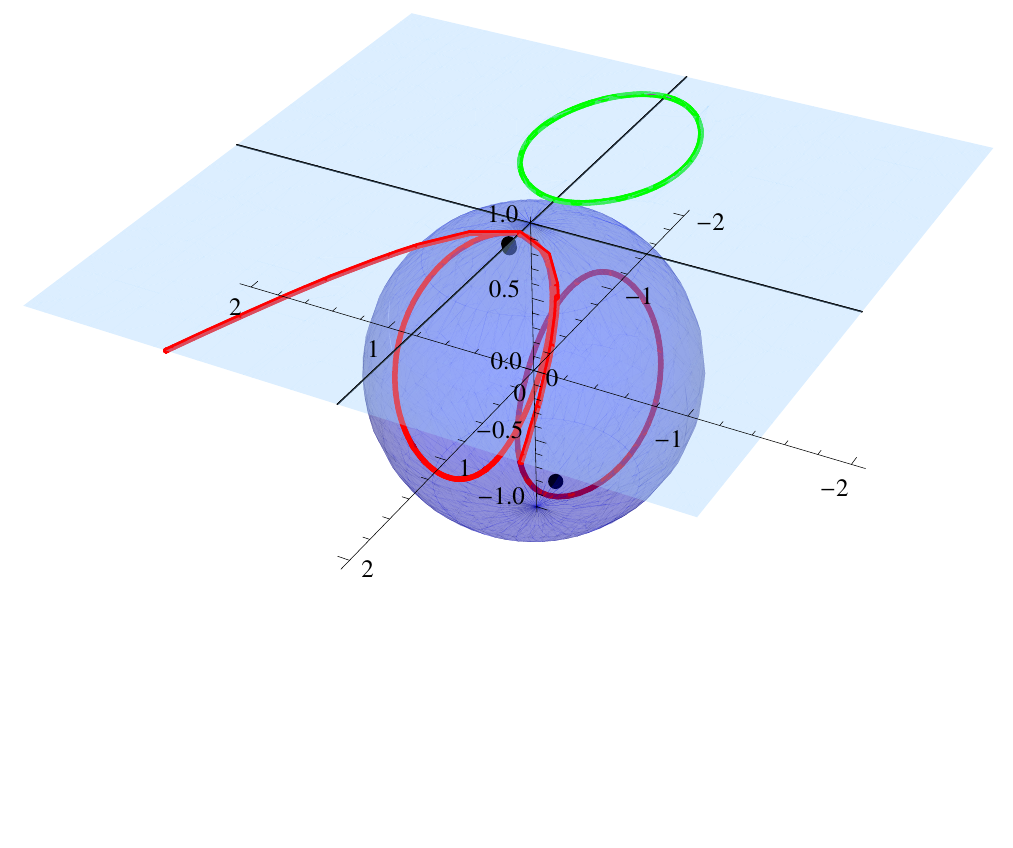} \includegraphics[scale=0.55]{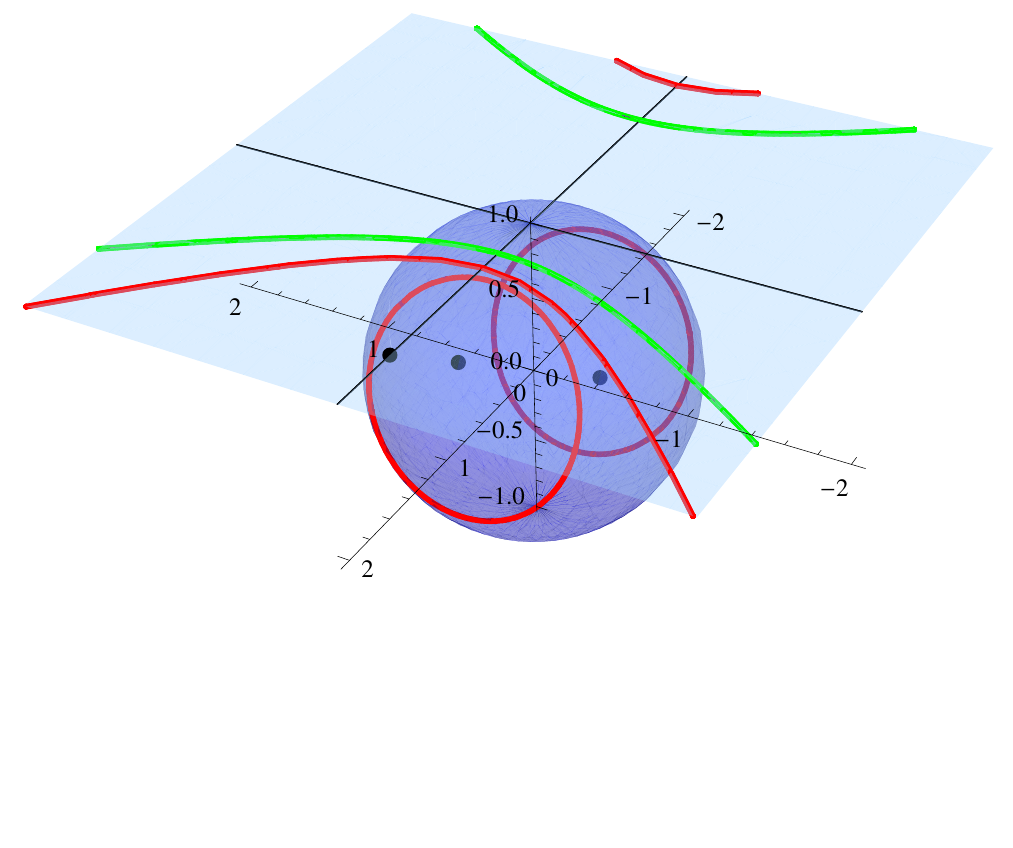}
\caption{The isoptic curve to elliptic parabola with parameters: \newline Parabola: $p=0.25$; $\alpha=\pi/3$~ ~ ~ ~ ~ ~ ~ ~$p=1.5$, $\alpha=2\pi/3$}
\end{figure}

We consider in this paper the hyperbolic conic sections
with proper foci, the problem is timely
for the "generalized" conis cection types in the hyperbolic
plain (see: \cite{M81}). Moreover, similar questions are interesting for other plane
geometries e.g. in the Minkowsky plane.


\end{document}